\documentclass[11pt]{article}
\usepackage{mathrsfs}
\usepackage{latexsym,lineno}
\usepackage{algorithm}
\usepackage{enumerate}
\usepackage{algorithmic}
\usepackage{epsfig}
\usepackage{color}
\usepackage{amsmath}\usepackage{fleqn}\usepackage{verbatim}\usepackage{epsf}
\usepackage{amsthm}\usepackage{graphicx, float}\usepackage{graphicx}
\usepackage{amsfonts}
\usepackage{amssymb}\usepackage{graphpap}
\usepackage{epic}\usepackage{curves}
\usepackage{tikz}
\usepackage{subfigure}
\usepackage{mathrsfs}
\usepackage{pstricks}
\usepackage{hyperref}

\newtheorem{fact}{Theorem}

\newtheorem{theorem}[fact]{Theorem}

\newtheorem{corollary}[fact]{Corollary}
\newtheorem{lemma}[fact]{Lemma}

\numberwithin{equation}{section}
\newtheorem{proposition}[fact]{Proposition}

\title{\bf Anti-Ramsey number of matchings in $r$-partite $r$-uniform  hypergraphs\thanks {Research was partially supported by the National Natural Science Foundation of China (Nos. 11871329, 11971298)}}
\author { Yisai Xue$^{1}$,\, Erfang  Shan$^{2}$, \, Liying Kang$^{1}$\thanks{\em  Email address: lykang@shu.edu.cn (L. Kang),  xys16720018@163.com (Y. Xue), efshan@shu.edu.cn (E. Shan)} \\
{\small $^{1}$Department of Mathematics, Shanghai University,
Shanghai 200444, P.R. China}\\
{\small$^{2}$School of Management, Shanghai University,
Shanghai 200444, P.R. China}}

\date{}
\begin{document}

\maketitle

\begin{abstract}

An edge-colored hypergraph is rainbow if all of its edges have different colors. Given two hypergraphs $\mathcal{H}$ and $\mathcal{G}$, the anti-Ramsey number $ar(\mathcal{G}, \mathcal{H})$ of $\mathcal{H}$  in $\mathcal{G}$ is  the maximum  number of colors in a coloring of the edges of $\mathcal{G}$ so that there does not exist a rainbow copy of $\mathcal{H}$. Li et al. determined the anti-Ramsey number of $k$-matchings in  complete bipartite graphs. Jin and Zang showed the uniqueness of the extremal coloring. In this paper,  as a generalization of these results,  we determine the anti-Ramsey number $ar_r(\mathcal{K}_{n_1,\ldots,n_r},M_k)$ of $k$-matchings in complete $r$-partite $r$-uniform  hypergraphs and show the uniqueness of the extremal coloring. Also, we show that $\mathcal{K}_{k-1,n_2,\ldots,n_r}$ is the unique  extremal hypergraph for Tur\'{a}n number $ex_r(\mathcal{K}_{n_1,\ldots,n_r},M_k)$ and  show that $ar_r(\mathcal{K}_{n_1,\ldots,n_r},$ $M_k)=ex_r(\mathcal{K}_{n_1,\ldots,n_r},M_{k-1})+1$, which gives a multi-partite version result of \"Ozkahya and Young's conjecture.

\bigskip

\medskip

\noindent{\bf AMS (2000) subject classification:}  05C35
\end{abstract}

\section{Introduction}
An   edge-colored graph $G$ is called \textit{rainbow} if every edge of $G$ receives a different color.
Given two graphs $H$ and $G$, $ar(H,G)$ is defined to be  the maximum number of colors in a coloring of the edges of $H$ that has no rainbow copy of $G$. The number $ar(H,G)$ is called the \textit{anti-Ramsey number} of $G$ in $H$. When $H=K_n$, $ar(K_n,G)$ is the anti-Ramsey number of $G$.
Let $ex(H, G)$ denote the maximum number of edges that a subgraph of $H$ can have with no  subgraph isomorphic to $G$.

The study of anti-Ramsey number began  by Erd\H{o}s et al. \cite{Erd1973Anti} in 1970s.
In the original work, they  conjectured that $ar\left(K_{n}, C_{k}\right)=\left(\frac{k-2}{2}+\frac{1}{k-1}\right) n+O(1)$, and proved the conjecture when  $k = 3$.
After that, Alon \cite{1983On} proved the conjecture for $k = 4$.
Jiang, Schiermeyer and West (unpublished manuscript) proved the conjecture for $k \leq 7$.
Finally, Montellano-Ballesteros and Neumann-Lara \cite{2005An} completely proved the conjecture in 2005.

For matchings, Schiermeyer \cite{2004Rainbow} used a counting technique to determine  $ar(K_n,kK_2)$ for all $k\geq 2$ and $n\geq 3k+3$.
After that, Fujita et al. \cite{2009A} solved this problem for $k\geq 2$ and $n\geq 2k+1$.
In 2009, Chen et al. \cite{CHEN20093370} extended Schiermeyer's result to all $k\geq 2$ and $n\geq 2k$ by using the Gallai-Edmonds structure theorem.

Taking complete bipartite graphs as the host graphs, Li et al \cite{LI20092575} determined  $ar(K_{n_1,n_2},kK_2)$ for all $k\geq 1$.
Denote by  $B_{n,m}$ the set of all the $m$-regular bipartite graphs of order $2n$.  Li and Xu \cite{LI20091525} showed that $ar(B_{n,m},kK_2)=m(k-2)+1$ for $k\geq 2$, $m\geq 3$ and $n>3k-1$.

A {\em hypergraph} $\mathcal{H}=(V(\mathcal{H}), E(\mathcal{H}))$ is a finite set $V(\mathcal{H})$ of elements, called {\em vertices}, together with a finite set $E(\mathcal{H})$ of subsets of $V(\mathcal{H})$, called {\em
hyperedges} or simply {\em edges}. 
The \textit{union} of hypergraphs $\mathcal{G}$ and $\mathcal{H}$ is the graph $\mathcal{G}\cup\mathcal{H}$
with vertex set $V(\mathcal{G})\cup V(\mathcal{H})$ and edge set $E(\mathcal{G})\cup E(\mathcal{H})$.
If each edge of $\mathcal{H}$ has exactly $r$ vertices,  $\mathcal{H}$ is called $r$-\textit{uniform}.
For a subset $V'$ of $V(\mathcal{H})$, denoted by $\mathcal{H}[V']$ the subhypergraph of $\mathcal{H}$ induced by $V'$.
For $v\in V(\mathcal{H})$, we use $\mathcal{H}-v$ to denote $\mathcal{H}[V(\mathcal{H})\backslash\{v\}]$. 
For an edge $e$ in $E(\mathcal{H})$, denoted by $\mathcal{H}-e$ the hypergraph obtained by deleting $e$ from $\mathcal{H}$.
For a vertex $v\in V(\mathcal{H})$, the {\em degree} $d_{\mathcal{H}}(v)$ is defined as the number of edges  of $\mathcal{H}$ containing $v$. 
A vertex of degree zero is called an \textit{isolated} vertex.
For $u, v \in V(\mathcal{H})$, we define $d_{\mathcal{H}}(u, v)$ to be the number of edges of $\mathcal{H}$ containing $\{u, v\}$, and we call this number the \textit{co-degree} of $\{u, v\}$.
For a hypergraph $\mathcal{H}$, we denote the number of edges in $\mathcal{H}$ by $e(\mathcal{H})$.
 A  \textit{complete $r$-uniform hypergraph} is a hypergraph whose edge
set consists of all $r$-subsets of the vertex set.
A \textit{matching} in a hypergraph is a set of edges  in which no two edges have a common vertex.
We call a matching with $k$ edges a \textit{$k$-matching}, denoted by $M_k$.
An edge-colored hypergraph is called \textit{rainbow hypergraph} if the all of its edges have different colors.
The \textit{representing hypergraph} of a hypergraph $\mathcal{H}$ with an edge coloring $c$ is a spanning subhypergraph of $\mathcal{H}$ obtained by taking one edge of each color of $c$.
For an edge set $E\subseteq E(\mathcal{H})$, let $c(E)$ denote the set of colors of edges in $E$. For simplicity,
when $E=\{e\}$ and $E=E(\mathcal{H})$, we use $c(e)$ and $c(\mathcal{H})$
  instead of $c(\{e\})$ and $c(E(\mathcal{H}))$, respectively.

Let $n_1, n_2, \ldots, n_r$ be integers and $V_1, V_2, \ldots, V_r$ be disjoint vertex sets with $|V_i|=n_i$  for each $i=1, 2, \ldots, r$.
A \textit{complete $r$-partite $r$-uniform hypergraph} on vertex classes $V_{1}, V_{2}, \ldots, V_{r}$, denoted by $\mathcal{K}_{n_1,\ldots,n_r}$,  is defined to be the $r$-uniform hypergraph
whose edge set consists of all the $r$-element subsets $S$ of $V_{1}\cup \cdots \cup V_{r}$ such that $|S\cap V_i|=1$ for all $i=1,2, \ldots,r$.

Given two hypergraphs $\mathcal{H}$ and $\mathcal{G}$,
the \textit{anti-Ramsey number} of $\mathcal{H}$  in $\mathcal{G}$,  denoted by $ar(\mathcal{G}, \mathcal{H})$, is  the maximum  number of colors in a coloring of the edges of $\mathcal{G}$ with no rainbow copy of $\mathcal{H}$. When $\mathcal{G}$ is an $r$-uniform complete hypergraph on $n$ vertices, $ar_r(\mathcal{G}, \mathcal{H})$ is  the  anti-Ramsey number of $\mathcal{H}$.
The \textit{Tur\'{a}n number} $ex_r(\mathcal{G},\mathcal{H})$ is the maximum number of edges in an $\mathcal{H}$-free subhypergraph of $\mathcal{G}$, where $\mathcal{H}$-free hypergraph is one which contains no $\mathcal{H}$ as a subhypergraph.

Gu et al. \cite{gu2020anti} determined the anti-Ramsey numbers
of linear paths/cycles and loose paths/cycles in hypergraphs for sufficiently large $n$ and gave bounds on the anti-Ramsey numbers of Berge paths/cycles.
For the anti-Ramsey number of matchings in hypergraphs, \"Ozkahya and Young \cite{OZKAHYA20132359}  stated a conjecture that $ar_r(\mathcal{K}_n,M_{k})=ex_r(\mathcal{K}_n,M_{k-1})+1$  for all $n>sk$ and proved the conjecture when $k=2,3$ and $n$ is sufficiently large.
Recently, Frankl and Kupavskii \cite{2019frankl} proved that the conjecture is true for $n\geq rk+(r-1)(k-1)$ and $k\geq 3$. Jin \cite{2021jin} determined the exact
value of the anti-Ramsey number of a $k$-matching in a complete tripartite $3$-uniform
hypergraph.

Take a subhypergraph $\mathcal{K}_{k-2, n_2,\ldots, n_r}$ of $\mathcal{K}_{n_1,n_2,\ldots,n_r}$. Color the edges of $\mathcal{K}_{k-2, n_2,\ldots, n_r}$ by distinct colors and color the remaining edges of $\mathcal{K}_{n_1,\ldots,n_r}$ by a new color. Denote by $\phi_r$ the obtained $((k-2)n_2\cdots n_r+1)$-edge-coloring of $\mathcal{K}_{n_1,\ldots,n_r}$.
Li, Tu and Jin \cite{LI20092575} determined the following results in complete bipartite graphs.

\begin{theorem}[\cite{LI20092575}]\label{li}
For $n_2\geq n_1\geq k \geq 1$,
$$ex(K_{n_1,n_2},kK_2)=(k-1)n_2.$$
Moreover, $K_{k-1,n_2}$ is the unique such extremal graph.
\end{theorem}

\begin{theorem}[\cite{LI20092575}]\label{thm3}
For $n_2\geq n_1\geq k\geq 3$,
$$ar(K_{n_1,n_2},kK_2)=(k-2)n_2+1.$$
\end{theorem}

In addition to the anti-Ramsey number, another interesting problem posed by Erd\H os  is the uniqueness of the extremal coloring.
In \cite{2017Anti}, Jin and Zang obtained the following result.

\begin{theorem}[\cite{2017Anti}]\label{Jin}
For $n_2\geq n_1\geq k\geq 3$, every $((k-2)n_2+1)$-edge-coloring except for $\phi_2$ of $K_{n_1,n_2}$ contains a rainbow $k K_{2}$.
\end{theorem}

The following proposition provides a lower and upper bound for $ar_r(\mathcal{K}_{n_1,\ldots,n_r},M_k)$, and the proof of which is similar to that of \cite{2004Rainbow}.

\begin{proposition}\label{prop}
$ex_r(\mathcal{K}_{n_1,\ldots,n_r},M_{k-1})+1 \leq ar_r(\mathcal{K}_{n_1,\ldots,n_r},M_k)\leq ex_r(\mathcal{K}_{n_1,\ldots,n_r},M_k).$
\end{proposition}

\noindent
\textbf{Proof.} The upper bound is clear.
For the lower bound, let $\mathcal{H}_{0}$ be an extremal hypergraph for $ex_r(\mathcal{K}_{n_1,\ldots,n_r},M_{k-1})$ and color all edges of $\mathcal{H}_{0}$ differently and all the edges in $E\left(\mathcal{K}_{n_{1}, \ldots, n_{r}}\right) \backslash E\left(\mathcal{H}_{0}\right)$ with one extra color.
The hypergraph $\mathcal{K}_{n_{1}, \ldots, n_{r}}$ with
this coloring does not contain a rainbow $k$-matching. The result follows. \qed

The proposition provides a lower bound for $ar_r(\mathcal{K}_{n_1,\ldots,n_r},M_k)$. In this paper we will show that $ar_r(\mathcal{K}_{n_1,\ldots,n_r},$ $M_k)=ex_r(\mathcal{K}_{n_1,\ldots,n_r},M_{k-1})+1$, which gives a multi-partite version result of \"Ozkahya and Young's conjecture.

In \cite{2020Tur}, Liu and Wang determined $ex_r(\mathcal{K}_{n_1,\ldots,n_r},M_k)$.

\begin{theorem}[\cite{2020Tur}]\label{ex}
For $n_r\geq n_{r-1}\geq \cdots\geq  n_1\geq k\geq  1$,
$$ex_r(\mathcal{K}_{n_1,\ldots,n_r},M_k)=(k-1)n_2\cdots n_r.$$
\end{theorem}

 We will show that $\mathcal{K}_{k-1,n_2,\ldots,n_r}$ is the unique  extremal hypergraph in Theorem \ref{ex}.

The following result is very useful for us.

\begin{theorem}\label{lemma}
For $n_r\geq n_{r-1}\geq \cdots\geq  n_1\geq k\geq  1$,
every subhypergraph  of $\mathcal{K}_{n_1,\ldots,n_r}$ with $(k-1)n_2\cdots n_r$ edges and without isolated vertices, except for $\mathcal{K}_{k-1,n_2,\ldots,n_r}$,  contains a $k$-matching.
\end{theorem}

Motivated by Theorem \ref{thm3}, one may naturally ask what is the maximum number of colors in a complete $r$-partite $r$-uniform hypergraph without a rainbow $k$-matching, where $r\geq 3$.
This paper  focus on  the anti-Ramsey number of $k$-matchings in complete $r$-partite $r$-uniform  hypergraphs. The following are our main results.

\begin{theorem}\label{thm1}
(i) For $n_r\geq n_{r-1}\geq \cdots \geq n_1\geq 3$,
$$ar_r(\mathcal{K}_{n_1,\ldots,n_r},M_2)=1.$$
(ii) For $n_1=2$,  let $t$ be the largest integer such that $n_t=n_1=2$. Then
$$ar_r(\mathcal{K}_{n_1,\ldots,n_r},M_2)=2^{t-1}.$$
\end{theorem}

\begin{theorem}\label{thm2}
For $n_r\geq n_{r-1}\geq \cdots\geq  n_1\geq 2k-1$ and $k\geq 3$,
$$ar_r(\mathcal{K}_{n_1,\ldots,n_r},M_k)=(k-2)n_2\cdots n_r+1.$$
Moreover, every $((k-2)n_2\cdots n_r+1)$-edge-coloring except for $\phi_r$ of $\mathcal{K}_{n_1,\ldots,n_r}$ contains a rainbow $k$-matching.
\end{theorem}

Combining Theorems \ref{ex}, \ref{thm1} and \ref{thm2}, we have the following corollary.

\begin{corollary}
For $n_r\geq n_{r-1}\geq \cdots\geq  n_1\geq 2k-1$ and $k\geq 2$,
$$ar_r(\mathcal{K}_{n_1,\ldots,n_r},M_k)=ex_r(\mathcal{K}_{n_1,\ldots,n_r},M_{k-1})+1.$$
\end{corollary}

\section{Proofs of Theorems \ref{lemma} and  \ref{thm1}}

\noindent\textbf{Proof of Theorem \ref{lemma}}.
We use induction on $r$.
The base case of $r=2$ is true for all $n_2\geq n_1\geq k$ by Theorem \ref{li}.
Suppose that the assertion holds for all $r'<r$.
Assume that $\mathcal{G}$ is a subhypergraph  of $\mathcal{K}_{n_1,\ldots,n_r}$ with $(k-1)n_2\cdots n_r$ edges and without isolated vertices, and does not contain a $k$-matching.
Let $V_s=\{v_{s1},v_{s2},\ldots,v_{sn_s}\}$ for $s=1,2,\ldots,r$. We consider two different cases.

\noindent\textbf{Case 1.} $n_1=n_2$.

For $1\leq i, j\leq n_1$, let $F_{i, j}=\big\{\{v_{1i},v_{2j},w_3,\ldots,w_{r}\}\in E(\mathcal{G})|~w_s\in V_{s} \text{ for } 3\leq s\leq r \big\}$ and $F_i=F_{i, 1}\cup F_{i+1,2}\cup\cdots\cup F_{i+n_1-1,n_1}$,  where $F_{i, j}=F_{i-n_1,j}$ if $i>n_1$.

For each $F_i$, $i=1,2,\ldots, n_1$, we construct an $(r-1)$-partite  $(r-1)$-uniform hypergraph $\mathcal{G}_{i}$  on vertex classes $V_1, V_3, \ldots, V_r$,
and $e=\{v_{1l}, w_3, \ldots, w_r \}$ is an edge of $\mathcal{G}_{i}$ if and only if $e'=\{v_{1l}, v_{2l'}, w_3, \ldots, w_r \}$ is an edge of $F_i$, where $l-l'\equiv i-1 \pmod {n_1}$.
Therefore, there is a bijection between $F_i$ and $E(\mathcal{G}_{i})$.
Note that if two edges $e_1$ and $e_2$ in $\mathcal{G}_{i}$ are independent,
then the corresponding edges $e_1'$ and $e_2'$ in $F_i$ are also independent. Then we have   the following fact.

\noindent\textbf{Fact A.} Any matching in $\mathcal{G}_{i}$ corresponds to a matching in $F_i\subseteq E(\mathcal{G})$.

First, we prove the following claims.

\noindent\textbf{Claim 1.} For $i\neq j$, $F_i\cap F_j=\emptyset$.

\noindent
\textbf{Proof.}
If there exists an edge $\{v_{1l},v_{2l'},w_3,\ldots,w_{r}\}\in F_i\cap F_j$, then $l-l'\equiv i-1 \pmod {n_1}$ and $l-l'\equiv j-1\pmod {n_1}$ which implies $i=j$.\qed

It follows from Claim 1 that  $e(\mathcal{G})=\sum\limits_{i=1}^{n_1}|F_i|=\sum\limits_{i=1}^{n_1}e(\mathcal{G}_{i})$.

\noindent\textbf{Claim 2.} For any $1\leq i\leq n_1$, $e(\mathcal{G}_{i})=(k-1)n_3\cdots n_r$.

\noindent
\textbf{Proof.}
First, we have $e(\mathcal{G}_{i})\leq (k-1)n_3\cdots n_r$. Otherwise, $\mathcal{G}_{i}$ contains a $k$-matching  by Theorem \ref{ex}, so  does $\mathcal{G}$
by Fact A, a contradiction.
Hence,
\begin{eqnarray*}
(k-1)n_1n_3\cdots n_r=e(\mathcal{G})=\sum\limits_{i=1}^{n_1}e(\mathcal{G}_{i})\leq n_1(k-1)n_3\cdots n_r,
\end{eqnarray*}
which implies that $e(\mathcal{G}_{i})= (k-1)n_3\cdots n_r$ for each $1\leq i\leq n_1$.\qed

According to Fact A and Claim 2,  $\mathcal{G}_{i}$ is a subhypergraph of $\mathcal{K}_{n_1, n_3, \ldots,n_r}$ with $(k-1)n_3\cdots n_r$ edges and does not contain a $k$-matching.
By the induction hypothesis,  $\mathcal{G}_{i}\cong \mathcal{K}_{k-1,n_3,\ldots,n_r}\cup (n_1-k+1)\mathcal{K}_1$ for $i=1,2,\ldots,n_1$.
Recall the construction of $\mathcal{G}_{i}$, we  deduce that $d_{\mathcal{G}}(v_{1i}, v_{2j})=0$ or $d_{\mathcal{G}}(v_{1i}, v_{2j})=n_3\cdots n_r$ for $1\leq i, j\leq n_1$.
Construct an auxiliary bipartite graph $G$ with bipartition $(V_{1},V_2)$, where $e_{ij}=v_{1i}v_{2j}\in E(G)$ if and only if $d_{\mathcal{G}}(v_{1i}, v_{2j})=n_3\cdots n_r$.
Then $e(\mathcal{G})=(k-1)n_2\cdots n_r$ implies that $e(G)=(k-1)n_2=(k-1)n_1$. 
We claim that there is no $k$-matching  in $G$. 
If there exists a $k$-matching $e_{i_1,j_1}$, $e_{i_2,j_2}$, $\ldots$, $e_{i_k,j_k}$ in $G$, 
we can find $k$ edges $e'_{i_1,j_1}$, $e'_{i_2,j_2}$, $\ldots$, $e'_{i_k,j_k}$ to form a $k$-matching in $\mathcal{G}$, where $e'_{i_l,j_l}=\{v_{1,i_l},v_{2,j_l},v_{3l},\ldots,v_{rl}\}$ for $l=1,2,\ldots,r$.
This contradicts the choice of $\mathcal{G}$.
It follows from Theorem \ref{li} that $G\cong K_{k-1,n_1}\cup (n_1-k+1)K_1$.
Without loss of generality, let $E(G)=\{v_{1i}v_{2j}|\ 1\leq i\leq k-1, 1\leq j\leq n_2\}.$
By the construction of $G$, every edge in $E(\mathcal{K}_{n_1,\ldots,n_r})$ containing $\{v_{1i},v_{2j}\}$ is an edge in $\mathcal{G}$ for $1\leq i\leq k-1$ and $1\leq j\leq n_2$.
Hence, $\mathcal{G}\cong\mathcal{K}_{k-1,n_2,n_3,\ldots,n_r}$ when $n_1=n_2$.

\noindent\textbf{Case 2.} $n_1<n_2$.

\noindent\textbf{Claim 3.} For $u\in V_2$, $d_{\mathcal{G}}(u)=(k-1)n_3\cdots n_r$.

\noindent
\textbf{Proof.}
If there exists a vertex $u\in V_2$ such that $d_{\mathcal{G}}(u)< (k-1)n_3\cdots n_{r}$, then
\begin{eqnarray*}
e(\mathcal{G}-u)=e(\mathcal{G})-d_{\mathcal{G}}(u)>(k-1)(n_2-1)n_3\cdots n_r.
\end{eqnarray*}
By Theorem \ref{ex}, $\mathcal{G}-u$ contains a $k$-matching, so does $\mathcal{G}$, a contradiction.
Hence, $d_{\mathcal{G}}(u)\geq (k-1)n_3\cdots n_{r}$ for all $u\in V_2$.
Note that
\begin{eqnarray*}
(k-1)n_2\cdots n_{r}=e(\mathcal{G})=\sum\limits_{u\in V_{2}}d_{\mathcal{G}}(u)\geq (k-1)n_{2}n_3\cdots n_{r}.
\end{eqnarray*}
We deduce that $d_{\mathcal{G}}(u)= (k-1)n_3\cdots n_{r}$ for all $u\in V_2$.\qed

Set $V_2'\subseteq V_2$ such that $|V_2'|=n_1$.
Let $\mathcal{G}'=\mathcal{G}[V_1,V_2',V_3,\ldots,V_r]$.
According to Claim 3, we have  $e(\mathcal{G}')=\sum\limits_{u\in V_{2}'}d_{\mathcal{G}}(u)=(k-1)n_1n_3\cdots n_r$.
It follows from Case 1 that
\begin{eqnarray*}
\mathcal{G}'\cong\mathcal{K}_{k-1,n_1,n_3,\ldots,n_r}\cup (n_1-k+1)\mathcal{K}_1.
\end{eqnarray*}
Combining this with  the arbitrariness of $V_2'$, we have
$d_{\mathcal{G}}(v_{1i}, v_{2j})=0$ or $d_{\mathcal{G}}(v_{1i}, v_{2j})=n_3\cdots n_r$ for $1\leq i\leq n_1, 1\leq j\leq n_2$.
Construct an auxiliary bipartite graph $G$ with bipartition $(V_{1},V_2)$, where $e_{ij}=v_{1i}v_{2j}\in E(G)$ if and only if $d_{\mathcal{G}}(v_{1i}, v_{2j})=n_3\cdots n_r$.
Then $e(\mathcal{G})=(k-1)n_2n_3\cdots n_r$ implies that $e(G) =(k-1)n_2$. 
We claim that there is no $k$-matching. Otherwise,
if there is a $k$-matching $e_{i_1,j_1},\ldots,e_{i_k,j_k}$ in $G$,  we can find $k$ edges in  $F_{i_1,j_1},\ldots,F_{i_k,j_k}$ to  form a $k$-matching in $\mathcal{G}$, a contradiction.
By Theorem \ref{li}, $G\cong K_{k-1,n_2}\cup (n_1-k+1)K_1$. Without loss of generality, let $E(G)=\{v_{1i}v_{2j}|\ 1\leq i\leq k-1, 1\leq j\leq n_2\}.$
By the construction of $G$, every edge in $E(\mathcal{K}_{n_1,\ldots,n_r})$ containing $\{v_{1i},v_{2j}\}$ is an edge in $\mathcal{G}$ for $1\leq i\leq k-1$ and $1\leq j\leq n_2$.
Hence, $\mathcal{G}\cong\mathcal{K}_{k-1,n_2,n_3,\ldots,n_r}$. \qed

We now turn to the proofs of the main results of this paper.
Theorem \ref{thm1} gives the  value of anti-Ramsey number of $k$-matching in complete $r$-partite $r$-uniform  hypergraphs when $k=2$.

\noindent\textbf{Proof of Theorem \ref{thm1}.}

\noindent\textbf{(i)} $n_1\geq 3$.

Suppose to the  contrary that  $\mathcal{H}$ is a complete $r$-partite $r$-uniform  hypergraph  colored by more than one color, and containing  no rainbow 2-matching.
Set $e,f\in E(\mathcal{H})$ such that $c(e)\neq c(f)$.
Clearly,  $E(\mathcal{H}-e)\cap E(\mathcal{H}-f)\neq \emptyset$ as $n_1\geq 3$.
Choose an edge $g\in E(\mathcal{H}-e)\cap E(\mathcal{H}-f)$, then $c(g)=c(e)$ as $\mathcal{H}$ does not contain a rainbow 2-matching. Similarly, we have $c(g)=c(f)$, then $c(e)=c(f)$,
contradicting the fact $c(e)\neq c(f)$.
Therefore, $|c(\mathcal{H})|=1$, i.e., $ar_r(\mathcal{K}_{n_1,\ldots,n_r},M_2)=1$.

\noindent\textbf{(ii)} $n_1=2$.

By the choice of $t$,
$n_1=n_2=\cdots=n_t=2$.
Let $\mathcal{H}\cong\mathcal{K}_{n_1,\ldots,n_r}$ be a complete $r$-partite $r$-uniform  hypergraph.
Let $R(\mathcal{H},t)=\{\{v_1,v_2,\ldots,v_t\}|~v_i\in V_i~\text{and}~i\in[t]\}$.
For any edge  $e=\{v_1,v_2,\ldots,v_r\}\in E(\mathcal{H})$, let $R(e,t)=\{v_1,v_2,\ldots,v_t\}$.
For $\alpha=\{v_1,v_2,\ldots,v_t\}\in R(\mathcal{H},t)$, let $\overline{\alpha}= \{\overline{v}_1, \overline{v}_2, \ldots, \overline{v}_t\}$, where $\overline{v}_i$ is the  remaining vertex in $V_i$ except $v_i$.
For $\alpha=\{v_1,v_2,\ldots,v_t\}\in R(\mathcal{H},t)$, let
\begin{eqnarray*}
E_{\{\alpha, \overline{\alpha}\}}=\{e\in E(\mathcal{H})\,|\,R(e, t)=\alpha\ \mbox{or}\   R(e, t)=\overline{\alpha}\}.
\end{eqnarray*}
Then we decompose the edge set of $\mathcal{H}$ into $E(\mathcal{H})=\cup_{\{\alpha, \overline{\alpha}\}\in Q}E_{\{\alpha, \overline{\alpha}\}},$ where $Q=\{\{\alpha, \overline{\alpha}\}|~ \alpha\in R(\mathcal{H}, t)\}$. It is easily checked that $|Q|=2^{t-1}$ and $E_{\{\alpha, \overline{\alpha}\}}\cap E_{\{\beta, \overline{\beta}\}}=\emptyset$ for any two different elements $\{\alpha, \overline{\alpha}\}, \{\beta, \overline{\beta}\}$ in $Q$.

We consider the following coloring of $\mathcal{H}$: For any two edges $e,f\in E(\mathcal{H})$, $c(e)=c(f)$ if and only if there exists an element $\{\alpha, \overline{\alpha}\}\in Q$ such that $e,f\in E_{\{\alpha, \overline{\alpha}\}}$. By the definition of the coloring of $\mathcal{H}$, any two independent edges have the same color and $|c(\mathcal{H})|=|Q|=2^{t-1}$. So $\mathcal{H}$ does not contain a rainbow $2$-matching.
Therefore, $ar_r(\mathcal{K}_{n_1,\ldots,n_r},M_2)\geq 2^{t-1}$.

Let $\mathcal{H}=\mathcal{K}_{n_1,\ldots,n_r}$ be an edge-colored complete $r$-partite $r$-uniform  hypergraph without rainbow 2-matchings.
To obtain the upper bound,  we consider the following  cases.

\noindent\textbf{Case 1.} $t=r$.

In this case, $E_{\{\alpha, \overline{\alpha}\}}=\{\alpha, \overline{\alpha}\}$. Then
$E(\mathcal{H})$ can be decomposed into $2^{r-1}$ pairs of independent edges, and each pair of edges must be colored by the same color.
Hence $|c(\mathcal{H})|\leq  2^{r-1}$, then $ar_r(\mathcal{K}_{2,\ldots,2},M_2)\leq 2^{r-1}$, which implies that $ar_r(\mathcal{K}_{2,\ldots,2},M_2) = 2^{r-1}$.

\noindent\textbf{Case 2.} $t<r$.

Recall that  $E(\mathcal{H})=\cup_{\{\alpha, \overline{\alpha}\}\in Q}E_{\{\alpha, \overline{\alpha}\}}$. We have the following claims.

\noindent\textbf{Claim 1.} For any $\{\alpha, \overline{\alpha}\}\in Q$, if $e,f\in E_{\{\alpha, \overline{\alpha}\}}$ and $R(e,t)=R(f,t)$, then $c(e)=c(f)$.

\noindent\textbf{Proof.} For any $\{\alpha, \overline{\alpha}\}\in Q$,
since $n_r\geq \cdots \geq n_{t+1}\geq 3$, there exists an edge $g= \{u_1,u_2,\ldots,u_r\}$ $\in E_{\{\alpha, \overline{\alpha}\}}$, such that $R(g,t)=\overline{R(e,t)}$ and $u_j\in V_j- e  \cup  f$ for $j=t+1,\ldots,r$.
Then
$g \cap e$ = $g \cap f=\emptyset$. Therefore, $c(e)=c(g)=c(f)$. Otherwise, $\mathcal{H}$ contains a rainbow 2-matching. \qed

\noindent\textbf{Claim 2.} For any $\{\alpha, \overline{\alpha}\}\in Q$, if $e,f\in E_{\{\alpha, \overline{\alpha}\}}$ and $R(e,t)=\overline{R(f,t)}$, then $c(e)=c(f)$.

\noindent\textbf{Proof.}
Since $n_r\geq \cdots \geq n_{t+1}\geq 3$, there exists an edge $g= \{u_1,u_2,\ldots,u_r\}\in E_{\{\alpha, \overline{\alpha}\}}$ with $R(g,t)=R(e,t)$ and $u_j\in V_j-f$ for $j=t+1,\ldots,r$. It follows from Claim 1 that $c(g)=c(e)$.
Note that $g \cap f=\emptyset$, we have $c(f)=c(g)$. Then $c(e)=c(f)$.  \qed

Combining these two claims, we conclude that $c(E_{\{\alpha, \overline{\alpha}\}})=1$ for any $\{\alpha, \overline{\alpha}\}\in Q$.
Therefore, $|c(\mathcal{H})|\leq2^{t-1}$, which implies that $ar_r(\mathcal{K}_{n_1,\ldots,n_r},M_2)\leq2^{t-1}$.
Thus $ar_r(\mathcal{K}_{n_1,\ldots,n_r},M_2)=2^{t-1}$.\qed

\section{Proof of Theorem \ref{thm2}}

In this section
we will determine the  value of the anti-Ramsey number of $k$-matchings in complete $r$-partite $r$-uniform  hypergraphs, and  also give the uniqueness of extremal coloring.

We need the following lemma.
\begin{lemma}\label{n1}
For $n_1\geq 2k-1$ and $k\geq 3$,
$$ar_r(\mathcal{K}_{n_1,\ldots,n_1},M_k)=(k-2)n_1^{r-1}+1,$$
and every $((k-2)n_1^{r-1}+1)$-edge-coloring except for $\phi_r$ of $\mathcal{K}_{n_1,\ldots,n_1}$ contains a rainbow $k $-matching.
\end{lemma}

\noindent
\textbf{Proof.}
We use induction on $r$.
The base case of $r=2$ is true for all $n_1\geq 2k-1$ by Theorem \ref{thm3} and Theorem \ref{Jin}.
Suppose that the lemma holds for all $r'<r$.
Assume, by way of contradiction, that $\mathcal{H}=\mathcal{K}_{n_1,\ldots,n_1}$ is a hypergraph with a $((k-2)n_1^{r-1}+2)$-edge-coloring and does not contain  a rainbow $k$-matching.

Let $V_s=\{v_{s1},v_{s2},\ldots,v_{sn_1}\}$,  $s=1,2,\ldots,r$. For $1\leq i, j\leq n_1$,
let
\begin{eqnarray*}
E_{i, j}=\big\{\{v_{1i},v_{2j},w_3,\ldots, w_{r}\}\in E(\mathcal{H})\,|\, w_s\in V_{s}\text{ for } 3\leq s\leq r \big\},
\end{eqnarray*}
and $E_i=E_{i, 1}\cup E_{i+1,2}\cup\cdots\cup E_{i+n_1-1,n_1}$, where $E_{i, j}=E_{i-n_1,j}$ if $i>n_1$.

For each $E_i, i=1,2,\ldots, n_1$, we construct a complete $(r-1)$-partite  $(r-1)$-uniform hypergraph $\mathcal{H}_{i}$ on vertex classes $V_1,V_3,\ldots,V_r$ such that $e=\{v_{1l}, w_3, \ldots, w_r \}$ is an edge of $\mathcal{H}_{i}$ if and only if $e'=\{v_{1l}, v_{2l'}, w_3, \ldots, w_r \}$ is an edge of $E_i$, where $l-l'\equiv i-1\pmod  {n_1}$, and we color $e$ by $c(e')$.
Therefore, there is a bijection between $E_i$ and $E(\mathcal{H}_{i})$ and $c(E_i)=c(\mathcal{H}_{i})$.
Note that if two edges $e_1$ and $e_2$ in $\mathcal{H}_{i}$ are independent,
then the corresponding edges $e_1'$ and $e_2'$ in $E_i$ are also independent. Then we have the following fact.

\noindent\textbf{Fact B.} Any rainbow matching in $\mathcal{H}_{i}$ corresponds to a rainbow matching in $E_i\subseteq E(\mathcal{H})$.


\noindent

 Obviously, $E(\mathcal{H})=\bigcup\limits_{i=1}^{n_1}E_i$. Then
\begin{eqnarray*}
\sum\limits_{i=1}^{n_1}|c(\mathcal{H}_{i})|=\sum\limits_{i=1}^{n_1}|c(E_i)|\geq |c(\mathcal{H})|=(k-2)n_1^{r-1}+2.
\end{eqnarray*}
Without loss of generality, we assume that $\mathcal{H}_{1}$ has the most colors in  $\mathcal{H}_{1}$, $\ldots$, $\mathcal{H}_{n_1}$.

\noindent\textbf{Claim 1.} $|c(\mathcal{H}_{1})|=(k-2)n_1^{r-2}+1$.

\noindent
\textbf{Proof.}
First, we have $|c(\mathcal{H}_{i})|\leq (k-2)n_1^{r-2}+1$ for $1\leq i\leq n_1$. Otherwise, by induction hypothesis,
$\mathcal{H}_{i}$ contains a rainbow $k$-matching, so does $\mathcal{H}$ by  Fact B, a contradiction.
If $|c(\mathcal{H}_{1})|\leq (k-2)n_1^{r-2}$, then
\begin{eqnarray*}
(k-2)n_1^{r-1}+2=|c(\mathcal{H})|\leq \sum\limits_{i=1}^{n_1}|c(\mathcal{H}_{i})|\leq n_1(k-2)n_1^{r-2},
\end{eqnarray*}
 a contradiction.
Hence, $|c(\mathcal{H}_{1})|=(k-2)n_1^{r-2}+1$.\qed

  We next show that there exists an integer $2\leq t\leq n_1$ such that $|c(\mathcal{H}_{t})|=(k-2)n_1^{r-2}+1$ and $c(\mathcal{H}_{1}) \cap c(\mathcal{H}_{t})=\emptyset$.
Otherwise,
\begin{eqnarray*}
|c(\mathcal{H})|\leq ((k-2)n_1^{r-2}+1)+(k-2)n_1^{r-2}\cdot(n_1-1)<(k-2)n_1^{r-1}+2,
\end{eqnarray*}
 contradicting the assumption of $\mathcal{H}$.
Since $\mathcal{H}_1$ and $\mathcal{H}_t$ do not contain a rainbow $k$-matching and with $(k-2)n_1^{r-2}+1$ colors, by
the induction hypothesis,  they are both colored by $\phi_{r-1}$.
Let $H_t$ be the representing hypergraph of $\mathcal{H}_t$.
Then
\begin{eqnarray*}
e(H_t)=|c(\mathcal{H}_{t})|>(k-2)n_1^{r-2},
\end{eqnarray*}
 so $H_t$ contains a $(k-1)$-matching by Theorem \ref{ex}.
Hence, there is a rainbow $(k-1)$-matching in $\mathcal{H}_{t}$, which corresponds to a rainbow $(k-1)$-matching, denoted by $M_{k-1}$, in $E_t$. The $(k-1)$-matching $M_{k-1}$ meets $k-1$ vertices in $V_1$ and $k-1$ vertices in $V_2$, respectively. 
As $n_1\geq 2k-1$, there exists an integer $s$, such that $v_{1s}\in V_1-V(M_{k-1})$ and $v_{2s}\in V_2-V(M_{k-1})$.
Then we can find an edge $e\in E_{s,s}\subseteq E_1$ such that $e\cap V(M_{k-1})=\emptyset$.
Recall that $c(\mathcal{H}_{1})\cap c(\mathcal{H}_{t})=\emptyset$, $c(\mathcal{H}_{1})=c(E_1)$ and $c(\mathcal{H}_{t})=c(E_t)$, then   $e\cup M_{k-1}$ is  a rainbow $k$-matching in $\mathcal{H}$, a contradiction.
Hence $ar_r(\mathcal{K}_{n_1,\ldots,n_1},M_k)\leq (k-2)n_1^{r-1}+1$.
Combining this with Proposition \ref{prop} and Theorem \ref{ex}, we have $ar_r(\mathcal{K}_{n_1,\ldots,n_1},M_k) = (k-2)n_1^{r-1}+1$.

 Now we prove the uniqueness of the extremal coloring. Suppose that  $|c(\mathcal{H})|=(k-2)n_1^{r-1}+1$  and $\mathcal{H}$ does not contain a rainbow $k$-matching.
Let $H$ be the representing hypergraph of $\mathcal{H}$.
Then $e(H)=(k-2)n_1^{r-1}+1$.
As the  discussion above, we construct $E_i\subseteq E(\mathcal{H})$ and a complete $(r-1)$-partite  $(r-1)$-uniform hypergraph $\mathcal{H}_{i}$ on vertex classes $V_1,V_3,\ldots,V_r$  for $i=1, 2, \ldots, n_1$. Without loss of generality, assume that $\mathcal{H}_{1}$ has the most colors in  $\mathcal{H}_{1}$, $\ldots$, $\mathcal{H}_{n_1}$.
Then $|c(\mathcal{H}_{1})|=(k-2)n_1^{r-2}+1$.

Notice that $\mathcal{H}_{1}$ does not contain a rainbow $k$-matching and
 $|c(\mathcal{H}_{1})|=(k-2)n_1^{r-2}+1$. By the induction hypothesis, $\mathcal{H}_{1}$ is colored by configuration $\phi_{r-1}$.
 Without loss of generality,  we assume that in $\mathcal{H}_{1}$,
 $E^\ast=\big\{\{v_{1i},w_3,\ldots,w_r\}\in E(\mathcal{H}_{1})~|~k-1\leq i\leq n_1 \text{ and } w_s\in V_s \text{ for } 3\leq s\leq r \big\}$
 are colored by one color and all the remaining edges of $\mathcal{H}_1$ colored by distinct colors.
Then, in $\mathcal{H}$, $c(E_{k-1,k-1})=\cdots=c(E_{n_1,n_1})=c(E^\ast)$, and the colors of the edges in $E_{1,1}\cup\cdots\cup E_{k-2,k-2}$ are different from each other.
Let $\mathcal{H}_0$ be the rainbow subhypergraph of $\mathcal{H}$ obtained by taking one edge of each color of $c(\mathcal{H})$ except $c(E^\ast)$ and such that $E_{1,1}\cup\cdots\cup E_{k-2,k-2}$ are contained in  $\mathcal{H}_0$.
Then $|E(\mathcal{H}_0)|=(k-2)n_1^{r-1}$. We have the following claim.

\noindent\textbf{Claim 2.} $\mathcal{H}_0\cong \mathcal{K}_{k-2,n_1,\ldots,n_1}$.

\noindent
\textbf{Proof.}
If $\mathcal{H}_0\ncong \mathcal{K}_{k-2,n_1,\ldots,n_1}$,
then there exists a $(k-1)$-matching, denoted by $M_{k-1}'$, in $\mathcal{H}_0$ by Theorem \ref{lemma}.
Notice that $M_{k-1}'$ meets at most $2(k-1)$ vertices in $V_1\cup V_{2}$ and $n_1\geq 2k-1$, we can find an edge $e_k'$ in $E_1$, such that $V(M_{k-1}')\cap e_{k}'=\emptyset$.
If $e_k'\in E_1\backslash (E_{1,1}\cup\cdots\cup E_{k-2,k-2})$, then
$c(M_{k-1}')\cap c(e_{k}')=\emptyset$.
If $e_k'\in E_{1,1}\cup\cdots\cup E_{k-2,k-2}$, then $e_k'\in E(\mathcal{H}_0)$.  We also 
 have $c(M_{k-1}')\cap c(e_{k}')=\emptyset$ as $\mathcal{H}_0$ is a rainbow subhypergraph of $\mathcal{H}$.
Then $e_k'\cup M_{k-1}'$ is a rainbow $k$-matching in $\mathcal{H}$, a contradiction.
\qed

By Claim 2,  there is a rainbow subhypergraph $\mathcal{H}_0\cong\mathcal{K}_{k-2, n_2,\ldots, n_r}$ of $\mathcal{H}$.
Combining the fact that all the edges in $E(\mathcal{H})\backslash E(\mathcal{H}_0)$ are colored by $c(E^*)$, we deduce that $\mathcal{H}$ is colored by configuration $\phi_r$.\qed
\vspace*{0.6cm}

We now prove Theorem \ref{thm2}.

\noindent\textbf{Proof of Theorem \ref{thm2}.}
We first prove that
$ar_r(\mathcal{K}_{n_1,\ldots,n_r},M_k)=(k-2)n_2\cdots n_r+1$ for  $n_1\geq 2k-1$ and $k\geq 3$.
Since $\phi_r$ is an edge coloring of $\mathcal{K}_{n_1,\ldots,n_r}$ such that there is no rainbow $k$-matching in $\mathcal{K}_{n_1,\ldots,n_r}$,
$ar_r(\mathcal{K}_{n_1,\ldots,n_r},M_k)\geq (k-2)n_2\cdots n_r+1$.

The upper bound is proven by induction on the number of  vertices of $\mathcal{K}_{n_1,\ldots,n_r}$.
The base case of $n_1=n_r$ is true by Lemma \ref{n1}, we now suppose that $n_1<n_r$.
Assume the assertion holds for all  $n_1',\ldots,n_r'$ such that $\sum_{i=1}^r n_i'<\sum_{i=1}^r n_i$.
Suppose that $\mathcal{H}=\mathcal{K}_{n_1,\ldots,n_r}$ is a hypergraph with a $((k-2)n_2\cdots n_r+2)$-edge-coloring, and containing no rainbow $k$-matchings.
Let $H$ be the representing hypergraph of $\mathcal{H}$. Then $e(H)=(k-2)n_2\cdots n_r+2$.
For any vertex $v\in V_{r}$, $H-v$ is  a subhypergraph of some representing hypergraph of $\mathcal{H}-v$.
Since $\mathcal{H}-v$ contains no rainbow $k$-matchings, by the induction hypothesis, we have
\begin{eqnarray*}
e(H-v)\leq |c(\mathcal{H}-v)|\leq (k-2)n_2\cdots n_{r-1}(n_{r}-1)+1.
\end{eqnarray*}
Then $d_H(v)=e(H)-e(H-v)\geq (k-2) n_2\cdots n_{r-1}+1$.
Thus,
\begin{eqnarray*}
e(H)&=&\sum_{v\in V_{r}}d_H(v)\geq n_{r}\cdot\big((k-2)n_2\cdots n_{r-1}+1\big)\\
&=&(k-2)n_2\cdots n_r+n_{r}>(k-2)n_2\cdots n_r+2,
\end{eqnarray*}
which is a contradiction.
Therefore, $ar_r(\mathcal{K}_{n_1,\ldots,n_r},M_k)=(k-2)n_2\cdots n_r+1$.

Next we prove the uniqueness of the extremal coloring.
Suppose that  $|c(\mathcal{H})|=(k-2)n_2\cdots n_r+1$ and $\mathcal{H}$ does not contain a rainbow $k$-matching.
Then $e(H)=(k-2)n_2\cdots n_r+1$. We have the following claim.

\noindent\textbf{Claim 1.} There exists a vertex $v^\ast\in V_r$ with $d_{H}(v^\ast)=(k-2)n_2\cdots n_{r-1}+1$, and  $d_{H}(v)=(k-2)n_2\cdots n_{r-1}$ for $v\in V_r\backslash\{v^\ast\}$.

\noindent
\textbf{Proof.}
For any vertex $v\in V_{r}$, $H-v$ is  a subhypergraph of some representing hypergraph of $\mathcal{H}-v$.
If $d_{H}(v)\leq (k-2)n_2\cdots n_{r-1}-1$, then
\begin{eqnarray*}
|c(\mathcal{H}-v)|\geq e(H-v)\geq (k-2)n_2\cdots n_{r-1}(n_r-1)+2=ar_r(\mathcal{K}_{n_1,\ldots,n_{r-1},n_{r}-1},M_k)+1.
\end{eqnarray*}
By the induction hypothesis, $\mathcal{H}-v$ contains a rainbow $k$-matching, a contradiction.
Thus,  $d_{H}(v)\geq (k-2)n_2\cdots n_{r-1}$ for all $v\in V_{r}$.
Note that
\begin{eqnarray*}
e(H)=\sum_{v\in V_{r}}d_{H}(v)= n_{r}\cdot(k-2)n_2\cdots n_{r-1}+1.
\end{eqnarray*}
Therefore, there is exactly one vertex, say $v^\ast$, with degree $(k-2)n_2\cdots n_{r-1}+1$, and the remaining vertices in $V_r$ with degree $(k-2)n_2\cdots n_{r-1}$.\qed

For  $v\in V_{r}\backslash\{v^\ast\}$, we have
\begin{eqnarray*}
e(H-v)=e(H)-d_H(v)=(k-2)n_2\cdots n_{r-1}(n_r-1)+1.
\end{eqnarray*}
Clearly, $H-v$ is  a subhypergraph of some representing hypergraph of $\mathcal{H}-v$.
Thus,
\begin{eqnarray*}
|c(\mathcal{H}-v)|\geq e(H-v)=(k-2)n_2\cdots n_{r-1}(n_r-1)+1.
\end{eqnarray*}
By the induction hypothesis, $|c(\mathcal{H}-v)|\leq (k-2)n_2\cdots n_{r-1}(n_r-1)+1$.
Hence,  we  deduce that $|c(\mathcal{H}-v)|=(k-2)n_2\cdots n_{r-1}(n_r-1)+1$.
By the induction hypothesis, $\mathcal{H}-v$ is colored by $\phi_r$ since $\mathcal{H}-v$ contains no rainbow $k$-matchings.
Due to the arbitrariness of $v$, for any two vertices $v_1,v_2\in V_{r}\backslash\{v^\ast\}$, $\mathcal{H}-v_1$ and $\mathcal{H}-v_2$ are colored by $\phi_r$.
Assume $S_1$ and $S_2$ are the $(k-2)$-subsets of $V_1$ such that the colors of the edges intersecting $S_1$ and $S_2$ are all different in   $\mathcal{H}-v_1$ and $\mathcal{H}-v_2$,  and edges not intersecting $S_1$ and $S_2$ are colored with a new color in $\mathcal{H}-v_1$ and $\mathcal{H}-v_2$, respectively. Then $S_1=S_2$. So all edges intersecting $S_1$ are colored with different colors  in $\mathcal{H}$ and  edges not intersecting $S_1$ are colored with a new color  in $\mathcal{H}$.
So
$\mathcal{H}$ is colored by $\phi_r$.
This completes the proof.
\qed


\end{document}